\documentclass[invmat,english,envcountsect]{svjour}
\usepackage[latin1]{inputenc}
\usepackage[T1]{fontenc}
\usepackage{amsmath,amsfonts}
\usepackage{graphics}
\usepackage{psfig,epsfig}
\def\R{\mathbb{R}}
\def\S{\mathbb{S}}
\def\N{\mathbb{N}}
\def\rn{\mathbb{R}^n}
\def\sn{\mathbb{S}^n}
\def\hn{\mathbb{H}^n}
\def\cotg{{\rm cotg}\,}
\def\sh{{\rm sh}\,}

\def\exp{{\rm exp}\,}
\def\Ric{{\rm Ric}\,}
\def\Vol{{\rm Vol}\,}
\def\Hess{{\rm Hess}\,}
\def\lap{\triangle\,}
\def\dis{\displaystyle}

\title{Eigenvalue pinching on convex domains in space forms}
\author{E. Aubry\thanks{Partially supported by FNRS Swiss Grant n°20-101469.}, J. Bertrand, B. Colbois}
\date{}
\institute{}
\DeclareMathOperator{\om}{\Omega}

\DeclareMathOperator{\ep}{\varepsilon}
\DeclareMathOperator{\vol}{ \mbox{ vol }}
\begin{document}
\maketitle
\abstract{In this paper, we show that the convex domains of $\mathbb{H}^n$ which are almost extremal for the Faber-Krahn or the Payne-Polya-Weinberger inequalities are close to geodesic balls. Our proof is also valid in other space forms and allows us to recover known results in $\rn$ and $\sn$.}

\section{Introduction}
\label{sec:intro}

This paper aims to study some optimal inequalities involving the first eigenvalues of the Dirichlet spectrum of convex domains in space forms, and to ask how stable they are. The paper essentially deals with the most intricate case of the hyperbolic space. 

The inequalities we are interested in are the Faber-Krahn inequality and the Payne-Polya-Weinberger inequality. The Faber-Krahn inequality asserts that among all bounded domains with the same volume in a given space form, the geodesic ball has the smallest first Dirichlet eigenvalue. Moreover, the geodesic ball is the unique minimizer (up to an isometry) among smooth domains. In this setting, such an inequality is stable if a bounded domain $\om$ whose $\lambda_1(\om)$ is close to $\lambda_1(B)$ ($B$ is a geodesic ball with the same volume as $\om$), is close for the Hausdorff distance to $B$ (up to an isometry). This general statement does not hold true, because it is possible to attach very long and thin tentacles to a ball without affecting significantly the volume and the spectrum. In fact, for Euclidean domains, weaker forms of stability have been established. One form is to prove that a domain whose first Dirichlet eigenvalue is close to the first eigenvalue of a suitable ball, ressembles a ball up to sets of small volume (see \cite{povel} for a precise statement). The other form is to consider only convex bodies, in this case the Faber-Krahn inequality is stable \cite{melas}. The stability of the Faber-Krahn inequality has also been established for convex domains in $\mathbb{H}^2$ and $\mathbb{S}^2$ \cite{avila}.

The first result of this paper is to prove the stability of the Faber-Krahn inequality for convex domains in a space form  of arbitrary dimension and arbitrary curvature. In the sequel, we will denote by
$X^1=(\sn,can)$, $X^0=(\rn,can)$ and $X^{-1}=(\hn,can)$ the space forms of curvature $1$, $0$ and $-1$
respectively. Except when stated otherwise, the results in this paper hold true for $\delta \in \{-1,0,1\}$. 

\begin{theorem}\label{stabFK}
Let $V_0>0$. Let $\lambda_1^*(V_0)$ be the first Dirichlet eigenvalue of a geodesic ball $B$ of volume $V_0$ in $X^\delta$.
 For any $\epsilon>0$, there exists $\eta>0$ such that, if $\Omega$ is a convex domain of volume $V_0$ in $X^\delta$ 
 and if $\lambda_1(\Omega)\leq\lambda_1^*(V_0)+\eta$ then, up to an isometry,
$$d_H(\om,B) \leq \epsilon,$$
 where $d_H$ denotes the Hausdorff distance. In the case $\delta=0$, we have $\eta=\eta'(\epsilon)V_0^{-2/n}$.
\end{theorem}

\begin{remark}We do not assume that the convex domains are bounded.
\end{remark}

The method developed is the same whatever the space form. Nevertheless, the case $\delta=-1$ is considerably harder. 
The primary difficulty is that the hyperbolic space comprises unbounded convex domains with finite volume therefore, have a discrete Dirichlet spectrum. This is contrary to the case of $\mathbb{R}^n$, where an upper bound of the  type ${\rm Diam}\Omega\leq C(\Vol\Omega,\lambda_1(\Omega),n)$ holds. To deal with this difficulty, we need to prove the thus for unsolved Faber-Krahn inequality for unbounded convex domains.

\begin{proposition}[Faber-Krahn Inequality]
Let $\om$ be a convex set in $X^\delta$ of finite volume $V_0$. The first Dirichlet eigenvalue of $\om$ satisfies
$$ \lambda_1(\om) \geq \lambda_1^*(V_0)$$
where $\lambda_1^*(V_0)$ denotes the first Dirichlet eigenvalue of a geodesic ball of volume $V_0$. Moreover, the equality $\lambda_1(\om)=\lambda_1^*(V_0)$  implies that $\om$ is isometric to a geodesic ball.
\end{proposition}
\begin{remark}The difficulty is in proving the case of equality.
\end{remark}

\smallskip

The second result of this paper concerns the stability of the Payne-Polya-Weinberger inequality (PPW inequality for short). This famous conjecture has been proved by M.S. Ashbaugh and R.D. Benguria \cite{ash1}.
\begin{theorem}[\cite{ash1}] Let $\om$ be a smooth bounded domain in Euclidean space (resp. a smooth domain included in a hemisphere in $\sn$). Then, the following inequality holds
$$ \frac{\lambda_2}{\lambda_1}(\om) \leq  \frac{\lambda_2}{\lambda_1}(B),$$
where $B$ is an arbitrary Euclidean ball (resp. a spherical ball such that $\Vol B=\Vol \Omega$). Moreover the equality is achieved if and only if $\om$ is isometric to a geodesic ball.
\end{theorem}
Let us notice that the ratio $\frac{\lambda_2(B)}{\lambda_1(B)}$ is scale-invariant in Euclidean space and that M. Ashbaugh and R. Benguria also showed in \cite{ash2}, that the ratio of the two first eigenvalues of a geodesic ball in $\sn$ is an increasing function of the radius $r$ (if $r\leq \pi/2$). Consequently, the PPW inequality follows directly from the theorem below.

\begin{theorem}[\cite{ash1,ash2,BL}]\label{genPPW}Let $\om$ be a smooth bounded domain in $X^{\delta}$
and such that $\om$ is included in an hemisphere if $\delta=1$. The second Dirichlet eigenvalue of $\om$ satisfies
$$ \lambda_2(\om) \leq \lambda_2(B)$$ 
where $B$ is a geodesic ball such that $\lambda_1(B)=\lambda_1(\om)$. Moreover, the equality holds if and only if $\om$ is isometric to $B$.
\end{theorem} 
 
It is shown in \cite{BL} that $\lambda_2/\lambda_1$ is a decreasing function of the radius of hyperbolic balls and that the PPW is false in $\hn$. This theorem can be seen as a generalized PPW inequality on space forms.

We prove the following stability results.

\begin{theorem}\label{stabB} Let $\om$ be a convex domain of $\rn$ or $\sn$, whose volume is equal to $V_0$. For any $\epsilon>0$ there exists $\eta>0$ such that for all $\om$ as above, the assumption $\frac{\lambda_2(\Omega)}{\lambda_1(\Omega)}\geq\frac{\lambda_2(B)}{\lambda_1(B)}-\eta$ implies
$$d_H(\om,B) \leq \ep,$$
where $B$ is a (well-centered) geodesic ball of volume $V_0$.
 \end{theorem}

\begin{remark}
The previous result was already known in Euclidean spa{\-}ce, it has been proved by A. Melas \cite{melas}.
\end{remark}
\begin{theorem}\label{stabBL}Let $\om$ be a convex domain of $X^{\delta}$ with $\lambda_1(\Omega)=\lambda$ ($\lambda > \frac{(n-1)^2}{4}$ if $\delta=-1$). For any $\ep>0$, there exits $\eta$ such that for all $\om$ as above, the assumption $ \lambda_2(\om) \geq \lambda_2^*(\lambda)-\eta$ implies
$$ d_H(\om,B)\leq \ep,$$
where $\lambda_2^*(\lambda)$ is the second Dirichlet eigenvalue of a (well-centered) geodesic ball $B$ of $X^\delta$ such that $\lambda_1(B)=\lambda$.
\end{theorem}

\begin{remark}We make no hypothesis about the volume of the convex domains we consider, not even the finiteness. This represents the main difference between this latter theorem and Theorems \ref{stabB} and \ref{stabFK}.
\end{remark}  

As for the Faber-Krahn inequality, it is necessary to generalize the PPW inequality to a more general setting, above all the characterization of the case of equality, in order to prove Theorems \ref{stabB} and \ref{stabBL} (see Theorems \ref{PPI} and \ref{conv2} for precise statements).

\medskip

The method developed to solve these stability questions is a rather general method, and based 
on the following abstract stability lemma.
The proof is straightforward, therefore omitted. 

Let $X$ be a topological space and $f : X\to\R$ be a function. 
We say that $f$ is {\sl coercive} if there exists a compact subset $K$ 
of $X$ such that $\inf_{X\setminus K}f>\inf_X f$ (we set $\inf_\emptyset f=+\infty$).

\begin{lemma}\label{lemmstab}
  Let $X$ be a topological space. If $f:X\to\R$ is coercive and lower semi-continuous 
  then $f$ is bounded below, reaches its minimal value and the set $M_f=f^{-1}\{\inf f\}$ of 
  its minima satisfies the following stability property:
for any neighborhood $U$ of $M_f$, there exists $\eta>0$ such that
$$f^{-1}\bigl(]-\infty,\inf f+\eta]\bigr)\subset U.$$
\end{lemma}

This lemma is very close to the so-called 
{\sl lower semi-continuity and compacity method}. This is 
typically used in calculus of variations to deal with the problem of minimisers existence (see \cite[Chapter 1]{Str}). 
It can be applied to a wide variety of problems (as large as the lower semi-continuity and compacity method). It does not, however, give an explicit $\eta$.

Our proof also shows that the infimum of the functional $\lambda_1$ (resp. $\frac{\lambda_1}{\lambda_2}$) on unbounded convex domains of $\hn$ with a given volume (resp. with a given $\lambda_1$) is strictly larger than those on bounded domains. To our knowledge, this is also a new result.

\bigskip
The paper is organized as follows:
\smallskip

In section \ref{sec:Dhaus} we define a metric on the space ${\cal C}$ of convex, bounded domains in $X^\delta$.

 In section \ref{sec:Lambdap} we show that the eigenvalues and volume functions are continuous on ${\cal C}$.

 In section \ref{Extr} we extend the classical Faber-Krahn and Payne-Polya-Weinberger (as its generalized version) inequalities to the set of convex unbounded domains. This level of generality is required in our proof even if this set is restricted to bounded convex domains for the proof of the coercivity.

Finally, we reduce the proof of the stability theorems to the proof of the coercivity of the functionals $\lambda_1$ and $\lambda_1/\lambda_2$ on the set of bounded convex domains of given volume (resp. given $\lambda_1$) and prove the coercivity of these functionals in sections \ref{Coer} and \ref{sec:Bengurian}. For that purpose, we prove several new qualitative results on the spectrum and the eigenfunctions of domains in space forms. For instance, we prove that a convex Euclidean domain with a spectral gap is bounded (hence has a discrete spectrum) and that its diameter is bounded from above by
$C(n)\bigl(\frac{1+\lambda_1}{\lambda_2-\lambda_1}\bigr)^{3/2}$ where $C(n)$ is a universal (explicitable) constant.

\section{A distance on convex domains}
\label{sec:Dhaus}

{\sl In the following, we set $s_1(t)=\sin t$, $s_0(t)=t$, $s_{-1}(t)=\sh t$ and $c_\delta=s'_\delta$. Let $x_0$ denote a fixed point in $X^\delta$.}

\begin{definition}
  Let ${\cal C}$ be the set of convex, bounded and open subsets $\Omega$ strictly included in $X^\delta$, which contain the point $x_0$.
\end{definition}

\begin{remark}
 The isometry group of $X^\delta$ acts transitively on $X^\delta$.
 \end{remark}
 \begin{remark}
Each proper, convex set of the sphere is included in a hemisphere. Hence, up to the sphere itself, all convex domains $\Omega$ in $\sn$ satisfy $vol(\om)\leq\Vol\sn/2$ and $\lambda_1(\om) \geq n$.
\end{remark}
So, Theorems \ref{stabFK}, \ref{stabB} and \ref{stabBL} are obvious in the case $\Omega=\sn$ and we just have to prove them for domains $\Omega\in{\cal C}$. 
\medskip

In the remaining part of this section we define a (proper) metric on ${\cal C}$. We chose to work with a metric which has a better behaviour than the usual Hausdorff metric with respect to the volume and the Dirichlet spectrum. To define our metric, we need some facts on support functions.

\subsection{Support functions}
\label{sec:Convbasic}

For any $\Omega\in{\cal C}$, the following function will be called {\sl the support of $\Omega$}.
$$
  \rho_\Omega:v\in S_{x_0}\mapsto\sup\{t\in\R_+/\,\exp_{x_0}(sv)\in\Omega\mbox{ for all }s\in[0,t]\}\in\R_+$$
where $S_{x_0}$ and ${\rm exp}_{x_0}$ are the set of unit tangent vectors and the exponential map of $X^\delta$ at $x_0$ respectively. Note that, even on $\sn$, we have $\rho_\Omega\leq R$ as soon as $\Omega\subset B(x_0,R)$.

 The properties of $\rho_\Omega$ needed subsequently are summarized in the following lemma.

\begin{lemma}\label{support}
 The function $\rho_\Omega$ is a Lipschitz function. Under the assumption $B(x_0,r)\subset\Omega\subset B(x_0,R)$, its Lipschitz constant is bounded above by $ s_\delta(R)\sqrt{\Bigl(\frac{s_\delta(R)}{s_\delta(r)}\Bigr)^2-1}$ if $\delta\neq1$, and by $\cotg r$ otherwise.\\
Moreover, the following equalities hold.\\
${~}\hskip1cm\overline{\Omega}={\rm exp}_{x_0}\{t.v/\,v\in S_{x_0},\,0\leq t\leq\rho_\Omega(v)\},$\\
${~}\hskip.6cm\quad\Omega={\rm exp}_{x_0}\{t.v/\,v\in S_{x_0},\,0\leq t<\rho_\Omega(v)\},$\\
${~}\hskip.6cm\quad\partial\Omega={\rm exp}_{x_0}\{\rho_\Omega(v).v/\,v\in S_{x_0}\}.$
\end{lemma}

\begin{proof}
Fix $y_0={\rm exp }_{x_0}(\rho_\Omega(u_0)u_0)\in \partial \om$ and consider the geodesic double cone centered at $y_0$ and tangent to the ball $B(x_0,r)$. We claim that for each $v\in S_{x_0}\setminus\{u_0\}$ close enough to $u_0$, the geodesic $\gamma_v(t)={\rm exp }_{x_0}(tv)$ meets the cone in exactly two points $Z(v),Z'(v)$. We also have  $l(d_{S_{x_0}}(v,u_0)) \leq \rho_\Omega (v) \leq L(d_{S_{x_0}}(v,u_0))$, where
$$l\bigl(d_{S_{x_0}}(v,u_0)\bigr)=\min \{d(x_0,Z(v)),d(x_0,Z'(v))\}$$
 and  
$$L\bigl(d_{S_{x_0}}(v,u_0)\bigr)=\max \{d(x_0,Z(v)),d(x_0,Z'(v))\}.$$

From elementary trigonometric computations (see appendix \ref{supp} for more details), we get
$$\displaylines{\liminf_{v\to u_0}\frac{\rho_\Omega(v)-\rho_\Omega(u_0)}{d_{S_{x_0}}(v,u_0)}\geq\liminf_{v\to u_0}\frac{l\bigl(d_{S_{x_0}}(v,u_0)\bigr)-l(0)}{d_{S_{x_0}}(v,u_0)}\hfill\cr
\hfill=l'(0)=-s_\delta\bigl(d(x_0,y_0)\bigr)\sqrt{\Bigl(\frac{s_\delta\bigl(d(x_0,y_0)\bigr)}{s_\delta(r)}\Bigr)^2-1}}$$
and
$$\displaylines{\limsup_{v\to u_0}\frac{\rho_\Omega(v)-\rho_\Omega(u_0)}{d_{S_{x_0}}(v,u_0)}\leq\limsup_{v\to u_0}\frac{L\bigl(d_{S_{x_0}}(v,u_0)\bigr)-L(0)}{d_{S_{x_0}}(v,u_0)}\hfill\cr
\hfill=L'(0)=s_\delta\bigl(d(x_0,y_0)\bigr)\sqrt{\Bigl(\frac{s_\delta\bigl(d(x_0,y_0)\bigr)}{s_\delta(r)}\Bigr)^2-1},}$$
which imply that $\rho_\Omega$ is Lipschitzian and give an upper bound on the Lipschitz constant, thanks to monotony properties of $s_\delta$.

The last three equalities of the statement follow easily from the continuity of $\rho_\Omega$ and standard properties of the exponential map.
\end{proof}

\subsection{A distance on convex bounded domains}

\begin{definition}
 We set $d(\Omega_1,\Omega_2)$ the metric defined on ${\cal C}$ by  $$d(\Omega_1,\Omega_2)=\|\ln\Bigl(\frac{\rho_{\Omega_1}}{\rho_{\Omega_2}}\Bigr)\|_\infty.$$
\end{definition}

\begin{proposition}\label{metrpropre}
  $\bigl({\cal C},d\bigr)$ is a proper metric space (i.e. every closed and bounded subset of $X$ is a compact set).
\end{proposition}

\begin{proof}
Let $(\Omega_i)_{i\in\N}$ be a bounded sequence in ${\cal C}$. Since there exist $r$ and $R$ such that $B(x_0,r)\subset\Omega_i\subset B(x_0,R)$ for every $i\in\N$, the functions $\rho_{\Omega_i}:S_{x_0}\to[r,R]$ are equicontinuous (thanks to Lemma \ref{support}). Consequently, the sequence $(\rho_{\Omega_i})_{i\in \mathbb{N}}$ converges uniformly on $S_{x_0}$ to a function $\rho_\infty$, up to an extraction. Since $r\leq\rho_\infty\leq R$, we also have that $\dis\lim_{i\to\infty}\|\ln\bigl(\rho_{\Omega_i}/\rho_\infty\bigr)\|_\infty=0$. We set $\Omega_\infty=\{\exp_{x_0}\bigl(t.v\bigr)/\,v\in S_{x_0},\, t\in[0,\rho_\infty(v)[\}$, which is a bounded, star-shaped domain. The equality $\rho_{\Omega_\infty}=\rho_\infty$ holds because $\rho_\infty$ is continuous and $\exp_{x_0}$  is a diffeomorphism of a neighbourhood of $B(0,R)$ onto a neighbourhood of $B(x_0,R)$. It remains to prove that $\Omega_\infty$ is convex. 

Let $y_1$ and $y_2$ be any pair of points in $\Omega_\infty$. There exists only one minimizing geodesic $\gamma$ from $y_1$ to $y_2$ in $X^\delta$ ($\Omega_\infty$ is an open set of a hemisphere in the case $\delta{=}1$). Since $y_1$ and $y_2$ are in $\Omega_j$ for all $j$ large enough, we easily infer that $\gamma\subset\overline{\Omega}_\infty=\{\exp_{x_0}\bigl(t.v\bigr)/\,v\in S_{x_0},\, t\in[0,\rho_\infty(v)]\}$.
So for any $r>0$ small enough, the union of the minimizing geodesic from $y_1$ (resp. from $y_2$) to a point of $B(y_2,r)$ (resp. of $B(y_1,r)$) is contained in $\overline{\Omega}_\infty$.
Since $y_1$ (resp. $y_2$) is in the injectivity domain of $y_2$ (resp. $y_1$), the union of this two sets is an open neighbourhood of $\gamma$ contained in $\overline{\Omega}_\infty$ and the result is proved.
 \end{proof}

\begin{corollary}\label{borne}
  For any $R\geq r>0$, the set of convex sets $\Omega$ in ${\cal C}$ such that $B(x_0,r)\subset\Omega\subset B(x_0,R)$ is a compact set.
\end{corollary}

\section{Continuity of the volume and the eigenvalues}
\label{sec:Lambdap}

 As proved in  \cite{GT}, any weak solution in $H^1_0(\om)$ of $\Delta u= \lambda u$ on a convex (in fact Lipschitzian) domain $\Omega$ belongs to $C^{\infty}(\om)\cap {\cal C}^0(\overline{\Omega})$ and is equal to $0$ on $\partial\Omega$. Moreover, the Dirichlet spectrum of  any open subset $\om$ of finite volume in $X^\delta$ is discrete \cite[Corollary 10.10]{shubin}. In this case, all the eigenvalues $(\lambda_k(\om)_{k\in\N^*}$ satisfy the min-max principle below.
$$\displaylines{\lambda_k(\Omega)= \inf\{m(E)/\,E \mbox{ subspace of }{\cal C}_c^\infty(\Omega), \dim E=k\}\cr
\mbox{where }m(E)=\sup_{f \in E}\frac{\int_{\om} |\nabla f|^2}{\int_{\om}f^2}.}$$
We will say that an arbitrary open set $\Omega$ {\sl has a spectral gap} if $\lambda_1(\Omega)<\lambda_2(\Omega)$ (where $\lambda_1(\om)$ and $\lambda_2(\om)$ are defined by the min-max principle). This implies that $\lambda_1(\Omega)$ is an eigenvalue of the Dirichlet problem and always occurs when the volume is finite.

\begin{proposition}\label{Lambdacont}
For any $k\geq1$, the following inequalities hold
$$\Bigl|\ln\Bigl(\frac{\lambda_k(\Omega_1)}{\lambda_k(\Omega_2)}\Bigr)\Bigr|\leq \Lambda_\delta\bigl[d(\Omega_1,\Omega_2),R\bigr]$$
and 
$$\Bigl|\ln\Bigl(\frac{\Vol\Omega_1}{\Vol\Omega_2}\Bigr)\Bigr|\leq \Lambda'_\delta\bigl[d(\Omega_1,\Omega_2),R\bigr]\text,$$
where 
$$\Omega_1\cup\Omega_2\subset B(x_0,R),\quad\Lambda_\delta(s,t)=\ln\Bigl[e^{2s}\Bigl(\frac{e^{2s}s_\delta(te^{-2s})}{s_\delta(t)}\Bigr)^{\delta(n-1)}\Bigr],$$
$$\Lambda'_1(s,t)=\Lambda'_0(s,t)=ns\quad\mbox{and}\quad\Lambda'_{-1}=\ln\Bigl[e^{ns}\Bigl(\frac{e^{-s}\sinh(t)}{\sinh(e^{-s}t)}\Bigr)^{n-1}\Bigr].$$
\end{proposition}

\begin{proof}
In the case $\delta=1$, we denote by $y_0$ the antipodal point of $x_0$ in $\sn$. For $\lambda\in]0,1]$, we define the  map 
\begin{eqnarray}\label{dilatations}
  H_\lambda:X^\delta\; (resp.\; X^1\setminus\{y_0\})&\to&X^\delta\\
\exp_{x_0}(tv)&\mapsto&\exp_{x_0}(\lambda tv) \nonumber
\end{eqnarray}
Set $d=d(\Omega_1,\Omega_2)$.
Since $H_{e^{-d}}(\Omega_1)\subset\Omega_2$ we just have to bound the quotient $\lambda_k\bigl(H_\lambda(\Omega_1)\bigr)/\lambda_k(\Omega_1)$ for $\lambda=e^{-d}$.

For that purpose, we define a linear injective map $\Phi_\lambda:{\cal C}_c^\infty(\Omega)\to{\cal C}_c^\infty\bigl(H_\lambda(\Omega)\bigr)$ by $\Phi_\lambda(f)=f\circ H_{1/\lambda}$. Easy computations involving Jacobi fields give
$$\lambda\Bigl(\inf_{t\in]0,R]}\frac{s_\delta(\lambda t)}{s_\delta(t)}\Bigr)^{n-1}\|f\|^{~}_1\leq\|\Phi_\lambda(f)\|^{~}_1\leq \lambda\Bigl(\sup_{t\in]0,R]}\frac{s_\delta(\lambda t)}{s_\delta(t)}\Bigr)^{n-1}\|f\|_1^{~}$$
$$\frac{\bigl|d\bigl(\Phi_\lambda(f)\bigr)\bigr|^2(x)}{\Phi_\lambda\bigl(|df|^2\bigr)(x)}\leq\max\Bigl(\frac{1}{\lambda^2},\frac{s^2_\delta(d(x_0,x)/\lambda)}{s_\delta^2(d(x_0,x))}\Bigr)$$
The first inequality applied to $f\equiv1$ gives the volume estimate. The two inequalities imply 
$$\frac{\|d\bigl(\Phi_\lambda(f)\bigr)\|_2^2}{\|\Phi_\lambda(f)\|_2^2}\leq e^{\Lambda_\delta(d,R)}\frac{\|df\|_2^2}{\|f\|_2^2}.$$
Using the min-max principle, we obtain
$$\lambda_k(\Omega_2)\leq\lambda_k\bigl(H_{e^{-d}}(\Omega_1)\bigr)\leq e^{\Lambda_\delta(d,R)}\lambda_k(\Omega_1)\text.$$
\end{proof}

\section{Extremal convex domains}\label{Extr}
\subsection{Schwarz symmetrization on noncompact domains}
The aim of this paragraph is to recall some basic properties of the Schwarz symmetrization. However we will not assume as usual, that the domain to be symmetrized is bounded. To replace this property, some additional assumptions on the functions to be symmetrized are sometimes needed.

\begin{definition}[Schwarz symmetrization] Let $f$ be a nonnegative function defined on an open set $\om$ in the space form $X^{\delta}$. Let $\mu_f$ be the distribution function defined for $s \geq 0$, by $\mu_f (s) = vol (\{f>s\})$ and let $V:\, r \mapsto Vol (B(r))$ ($r \geq 0$). The nonincreasing Schwarz symmetrization of $f$ is 
$$ f^* = \mu_f^{\sharp} \circ V \circ d_{x_0},$$
where $d_{x_0}(x)=d(x_0,x)$ and $.^{\sharp}$ refers to the right inverse function of a nonincreasing function (i.e. $u^\#(s)=\inf\{t\geq 0/ u(t)\leq s\}$). If the volume of $\om$ is finite, the Schwarz nondecreasing symmetrization of $f$ is defined by 
$$ f_* = \mu_f^{\sharp} \circ H \circ d_{x_0},$$
where $H : r \mapsto \vol(\om)-V(r)$.

These symmetrized functions 
satisfy 
\begin{equation}\label{equi}
\mu_{f^*}=\mu_{f_*}=\mu_f.
\end{equation}
\end{definition}    
\begin{remark}     
For more details on symmetrization, we refer to \cite{chong,kawohl,berard}.
\end{remark}
\begin{proposition}\label{symme}Let $\om$ be an open set of finite volume in the space form $X^\delta$.

If $u$ is in $L^2(\om)$ then $u^*$ is in $L^2(\om^*)$ and
\begin{equation}\label{mapl2}
 ||u||_{L^2(\om)}=||u^*||_{L^2(\om^*)}.
\end{equation}
In addition, the following inequalities hold
\begin{equation}\label{inel2}
 \int_{\om^*} f_*g^* \leq \int_{\om} fg \leq  \int_{\om^*} f^*g^* 
\end{equation}
for every nonnegative measurable functions $f,g$ on $\om$.

If $u$ is now in $H^1_0(\om)$ then $u^* \in H^1_0(\om^*)$ and
\begin{equation}\label{maph1}
 \int_{\om^*} |\nabla u^*|^2 \leq \int_{\om} |\nabla u|^2.
\end{equation} 
\end{proposition}
\begin{proof}
The proof of the statement (\ref{mapl2}) is an immediate consequence of (\ref{equi}), the inequality (\ref{inel2}) is easy to check for simple functions and the general case follows by density \cite{kawohl}. The proof of (\ref{maph1}) also relies on an approximation argument, a suitable dense subset is introduced in the lemma below. The assumption on the volume is then used to conclude, using Rellich's Theorem on the symmetrized ball and the following inequality which is a direct consequence of (\ref{mapl2}) and (\ref{inel2}).
$$ || u^*-v^*||_{L^2(\om^*)} \leq  || u-v||_{L^2(\om)}.$$
\end{proof} 
 
\begin{lemma}\label{lemtec}
Let $f$ be a smooth nonnegative function in $H^1_0(\om)$, which is zero on $\partial\Omega$ and in ${\cal C}^0(\overline{\Omega})$, where $\om$ is an open set of finite volume in $X^{\delta}$. Suppose that the level sets of $f$ are compact sets (except maybe $\{f=0\}$) of measure zero. Under these assumptions,  $\mu_f^{\sharp}$ is absolutely continous, the symmetrized function $f^*$ is in $H^1_0(\om^*)$ and satisfies (\ref{maph1}).

Moreover, in case of equality in (\ref{maph1}), the open set $\{f>0\}$ is a ball.
\end{lemma}

\begin{proof}
Let $Reg(f)$ be the set of regular points of $f$ which are included in $\{x \in \om; f>0\}$. By assumption, $Reg(f)$ is an open set of full measure in $\{x \in \om; f>0\}$. As a consequence, we deduce that $f^*$ is continuously differentiable on an open set of full measure of $\{f^*>0\}$ and satisfies the inequality (\ref{maph1}) thanks to the coarea formula and the isoperimetric inequality (we refer to \cite{BM} for more details). We conclude that $\{f>0\}$ is a ball using  a decreasing sequence of regular values which goes to 0 and the case of equality in the isoperimetric inequality. 
\end{proof}

\begin{remark} The set of functions which satisfy the assumptions of the lemma above contains the smooth functions with compact support and only nondegenerate critical points, therefore it is dense in $H^1_0(\om)$  (see \cite{BM} and references herein).
\end{remark} 
        
\begin{remark}\label{Sch} In the sequel, we will use the Schwarz symmetrization on convex domains of $\mathbb{H}^n$ whose the volume is not assumed to be finite. A priori, the nondecreasing Schwarz symmetrization cannot be defined in this setting, however the inequality 
$$ \int_{\om^*} f_*g^* \leq \int_{\om} fg $$
remains true for a function $f= F\circ d_{x_0}$, where $F$ is a nonnegative and nondecreasing bounded function such that $F$ is constant outside a compact set, if we define $f_*$ as 
$ 
f_*(x) = \left\{ \begin{array}{ll} (f|_{\om \cap B(x_0,r)})_* & \mathrm{ if}\ |x| < r \\
        ||f||_{\infty} & \mathrm{ otherwise}
\end{array}
\right.
$
for $r$ large enough.
\end{remark}

\subsection{Faber-Krahn Inequality} 
\label{subsec:spec}
In this section, we extend the Faber-Krahn inequality from the setting of smooth bounded domains to the setting of convex sets of finite volume. The main interest of the result below is the characterization of the case of equality.

\begin{proposition}[Faber-Krahn Inequality]\label{FK}
Let $\om$ be a convex set in $X^\delta$ of finite volume $V_0$. The first Dirichlet eigenvalue of $\om$ satisfies
$$ \lambda_1(\om) \geq \lambda_1^*(V_0)$$
where $\lambda_1^*(V_0)$ denotes the first Dirichlet eigenvalue of a geodesic ball with volume $V_0$. Moreover, the equality $\lambda_1(\om)=\lambda_1^*(V_0)$  implies that $\om$ is isometric to a geodesic ball.
\end{proposition}

\begin{remark}The inequality can be proved using smooth approximations of $\om$. However, the characterization of the case of equality without assuming that $\om$ is bounded, is crucial in our proof of Theorem \ref{stabFK}, when $\delta=-1$. Even when the domain is bound\-ed, some regularity on the boundary is needed to deduce the case of equality. Indeed each ball with closed sets of capacity zero removed, satisfies the case of equality.  
\end{remark}

\begin{proof}
The proof of the inequality  follows from Proposition \ref{symme} and does not rely on the convexity of $\om$.
As the volume of $\Omega$ is assumed to be finite, the Dirichlet spectrum of $\om$ is discrete \cite[Corollary 10.10]{shubin} and the eigenfunctions belong to $C^{\infty}(\om)\cap{\cal C}^0(\overline{\Omega})$ \cite[Corollary 8.11 and theorem 8.29]{GT}. To prove the case of equality, it is sufficient to prove that the first eigenfunction (denoted by $f_1$) satisfies the assumptions of Lemma \ref{lemtec}, which is a consequence of the lemma below. Indeed, thanks to this lemma and Sard's Theorem, the set of singular values of $f_1$ is a closed set of measure zero. Then, thanks to the fact that the function $\Delta f_1=\lambda_1 f_1$ is positive on $\om$, we deduce that each level set of the first eigenfunction is of measure zero.
\end{proof} 
\begin{remark}\label{rma}Let us remark that the assumption on the finiteness of the volume is used only to prove that the bottom of the spectrum is an eigenvalue. It is also true for the lemma below; we will use this fact in Paragraph \ref{sec:Ben}. 
\end{remark}

\begin{lemma}\label{KOMPAQ} Under the assumptions of Proposition \ref{FK}, the first Diri{\-}chlet eigenfunction $f_1$ on $\om$ can be assumed to be positive and proper: for all $s>0$, the set $f_1^{-1}([s,+ \infty[)$ is a compact set.
\end{lemma}

\begin{proof}
By the maximum principle, we can suppose $f_1$ to be positive. To prove the second assertion, set $y_0$ be a fixed point of $\Omega$, $R\geq1$ and $x_0\in\Omega\setminus B(y_0,2R)$. Recall (see for instance \cite{Cha}) that there exists a constant $C(n)$ such that,
$$\forall x_0\in X^\delta,\quad \forall v\in H^1_0\bigl(B(x_0,1)\bigr),\quad\|v\|_\frac{2n}{n-2}^2\leq C(n)\|dv\|_2^2\text.\quad(*)$$
Note that in dimension $n=2$, this inequality has to be replaced by $\|v\|_{4}^2\leq C\|dv\|_2^2$ in what follows. A standard Moser's iteration gives
\begin{equation}\label{itera}
  f_1^2(x_0)\leq A(n)\bigl(1+\lambda_1\bigr)^{\gamma(n)}\int_{B(x_0,1)} f_1^2
\end{equation}

(where $A(n)$ and $\gamma(n)$ are constants that depend only on the dimension $n$, see appendix \ref{iter} for a proof),
and from which we infer that
$$\sup_{\Omega\setminus B(y_0,2R)} f_1^2\leq A(n)\bigl(1+\lambda_1\bigr)^{\gamma(n)}\int_{\Omega\setminus B(y_0,R)} f_1^2.$$
This gives the compactness of the sets $f_1^{-1}\bigl([s,+\infty[\bigr)$ for all $s>0$ since $\int_{\Omega\setminus B(y_0,R)} f_1^2\to 0$ when $R\to\infty$ and $f_1$ is continuous on the convex set $\overline{\Omega}$ and is equal to $0$ on $\partial\Omega$.
\end{proof}

\subsection{Payne-Polya-Weinberger Inequality}\label{sec:Ben}
M.S. Ashbaugh and R.D. Benguria proved the Payne-Polya-Weinber{\-}ger conjecture for smooth bounded sets of Euclidean space and smooth sets included in a hemisphere of the sphere \cite{ash1,ash2}. We need to extend this inequality to (possibly non smooth) convex sets in the space form $X^\delta$ ($\delta \in \{0,1\}$). 
\begin{theorem}[Payne-Polya-Weinberger Inequality]\label{PPI}Let $\om$ be a convex set of finite volume $V_0$ in $X^{\delta}$ ($\delta \in \{0,1\}$). Under these assumptions, the following inequality is satisfied, 
$$ \frac{\lambda_2}{\lambda_1}(\om) \leq \frac{\lambda_2^*}{\lambda_1^*}(V_0).$$
Moreover, the equality is achieved if and only if $\om$ is isometric to a geodesic ball.
\end{theorem}

\begin{remark}\label{REM}
Actually, as in \cite{ash1,ash2}, the monotony properties of the ratio $\lambda_1(B)/\lambda_2(B)$ with respect to the radius of the geodesic ball $B$ of $X^\delta$ make this theorem a direct corollary of the following result.
\end{remark}

\begin{theorem}\label{conv2}  
Let $\om$ be a convex set in $X^\delta$ such that $\om \neq X^\delta$. Then the spectral gap of $\Omega$ is smaller or equal to $\lambda_2(B)-\lambda_1(B)$ (where $B$ is a geodesic ball such that $\lambda_1(B)=\lambda_1(\om)$). If the spectral gap is equal to $\lambda_2(B)-\lambda_1(B)$, $\om$ is isometric to a geodesic ball.
\end{theorem} 

Let us remark that contrary to the cases $\delta \in \{0,1\}$, the assumptions in Theorem \ref{conv2} do not imply an upper bound on the volume of $\om$ in $\hn$.

We will prove Theorems \ref{PPI} and \ref{conv2} simultaneously. The scheme of the proof is the same as in \cite{ash1,ash2,BL}, so we will mainly focuse on the extra arguments needed in our setting. The fisrt step of the proof is the following proposition.

\begin{proposition}\label{specgap}
  Let $\Omega$ be an open subset of $X^\delta$ (included in a hemisphere if $\delta=1$) with a spectral gap, $u_1$ an eigenfunction of $\Omega$ for the first eigenvalue and $g$ be a positive, piecewise ${\cal C}^1$ function on $[0,\infty[$ (and with $\liminf_{+\infty}g>0$ if $\Omega$ is not bounded). Then, there exists a point $x_m\in X^\delta$ such that
$$\lambda_2(\Omega)-\lambda_1(\Omega)\leq\frac{\int_\Omega b\bigl(d(x_m,y)\bigr)u_1^2(y)\,dy}{\int_\Omega g^2\bigl(d(x_m,y)\bigr)u_1^2(y)\,dy}$$
where $b=g'^2+\frac{n-1}{s_\delta^2}g^2$.
\end{proposition}

Note that for the proof of Theorem \ref{conv2}, we can suppose that the spectral gap is non zero.

\begin{proof}
The min-max principle implies that
$$ \lambda_2(\om)-\lambda_1(\om) \leq \frac{\int_{\om}|\nabla P|^2 u_1^2}{\int_{\om}P^2u_1^2},$$
for every non-zero function $P$ such that $Pu_1$ is in $H^1_0(\om)$ and \\$\int_{\om}Pu_1^2=0$.

The next step consists in choosing $n$ suitable test functions. For that purpose, we need the following lemma which extends a result of \cite{ash1,ash2,BL} (the proof is postponed to appendix \ref{proofcenter}).
\begin{lemma}\label{masscenter}
For any $u\in L^2(X^\delta)$ (with support in a hemisphere if $\delta=1$) and any $g:\R^+\to\R^+$ continuous (bounded and with $\liminf_{+\infty}g>0$ if $u$ has not compact support), there is $x\in X^\delta$ such that
$$\int_{X^\delta}g\bigl(d(x,y)\bigr)\frac{{\rm exp}_x^{-1}(y)}{d(x,y)}u^2(y)\,dy=0_{T_xX^\delta}$$
\end{lemma}

In order to construct test functions, we apply this lemma to $u=u_1.1_\Omega$ and $g$ a nonnegative, increasing and bounded function ($g$ will be specified later). The functions $P_i= g(r)X_i$, where $(r,X_i)$ are the geodesic coordinates at the point $x_m$ given by Lemma \ref{masscenter}, satisfy $\int_{\om}P_iu_1^2=0$ for every $i$. To conclude the proof of Proposition \ref{specgap}, we just have to sum the $n$ inequations given by the min-max principle applied to the $P_i$ ,and note that $\sum_i P_i^2=g^2$ and $\sum_i|\nabla P_i|^2=b$.

Now, we choose $g$ a radial function  such that the equality below holds.
$$ \lambda_2(B)-\lambda_1(B)= \frac{\int_Bb z^2}{\int_Bg^2z^2},$$
where $z$ is a positive first eigenfunction of $B$. It is shown in \cite{ash1,ash2,BL} that we can choose $g$ positive, nondecreasing and constant outside $B$ and such that $b$ is radial, positive and nonincreasing.
We recall that $B$ is such that $\lambda_1(B)=\lambda_1(\om)$. It remains to compare the spectral gaps. For that purpose, we first use properties of the Schwarz symmetrization (Proposition \ref{symme}). We get
$$ \int_{\om}b \,u_1^2 \leq \int_{\om^*}b^*u_1^{*2} \leq \int_{\om^*} b\,u_1^{*2}$$
$$ \int_{\om}g \,u_1^2 \geq \int_{\om^*}g_*u_1^{*2} \geq \int_{\om^*} g\,u_1^{*2}.$$
The inequality involving the nonincreasing Schwarz symmetrization is valid without any assumption on the volume, thanks to remark \ref{Sch}. We conclude using the Chiti comparison result, which allows to compare $z$ with $u_1^*$ on $B$. This comparison result is valid as soon as the first eigenfunction $u_1$ satisfies the assumptions of Lemma \ref{lemtec} (this has been established in the proof of Proposition \ref{FK} and does not rely on any assumption on the volume), we refer to \cite[pages 21-24]{ash1} for more details. Using the Chiti comparison result, we get \cite[page 607]{ash1}
$$\int_{\om^*} g \,u_1^{*2} \geq \int_B g\, z^2 \mbox{ and } \int_{\om^*} b\, u_1^{*2} \leq \int_B b\, z^2 $$
and this concludes the proof of the inequality. The case of equality follows from the characterization of the equality in the Chiti comparison result.
\end{proof}

\section{Coercivity of the functional $\lambda_1$}\label{Coer}

We show that $\lambda_1$ is coercive on appropriate subsets of ${\cal C}$. We first need a control of the in-radii of elements of ${\cal C}$.

\subsection{In-radius estimate in ${\cal C}$}
\label{sec:inrad}

For any bounded domain $\Omega$ in $X^\delta$, let ${\rm Inrad}\,(\Omega)$ be the maximum radius of a geodesic ball included in $\overline{\Omega}$ .

\begin{proposition}\label{inrad}
  Let $\Omega$ be a bounded convex set in $X^\delta$. Then
$${\rm Inrad}\,(\Omega)\geq\frac{\pi}{2\sqrt{\lambda_1(\Omega)+(n-1)}}\text.$$
\end{proposition}

This proposition has been proved by P.~Li and S.T.~Yau  \cite{LY} for smooth domains of nonnegative mean curvature (see appendix \ref{LY} for a proof). It can be readily extended to any (non smooth) convex domains: indeed, for any $\epsilon>0$ small enough, there exists a smooth convex domain $\Omega_\epsilon$ such that $H_{1-\epsilon}(\Omega)\subset\Omega_\epsilon\subset H_{1+\epsilon}(\Omega)$, where $H$ is the map defined by \eqref{dilatations}, p. \pageref{dilatations} (see \cite[Lemma 2.3.2]{Hor} for the Euclidean case and use the Klein projective model of the hyperbolic space and the open hemisphere, to infer this property in $\hn$ and $\sn$). The continuity of $\lambda_1$ on $\mathcal{C}$ allows to conclude.

\subsection{Coercivity of $\lambda_1$}

Subsequently, we denote by ${\cal C}_{V_0}$ the set of convex bounded domains $\Omega$ of $X^\delta$ with $\Vol \Omega=V_0$ and $B(x_0,{\rm Inrad}\,(\Omega))\subset\overline{\Omega}$ (note that ${\cal C}_{V_0}$ contains, up to isometry, all convex bounded domains of $X^\delta$ with volume $V_0$).

Combining Corollary \ref{borne} and Proposition \ref{inrad}, we get

 \begin{corollary}\label{bornes}
   For any $M>0$, the set of elements $\Omega$ of ${\cal C}$ (resp. ${\cal C}_{V_0}$) with $\Omega\subset B(x_0,M)$ and $\lambda_1(\Omega)\leq M$ is compact.
 \end{corollary}

\subsubsection{case $\delta{=}1$.}
\label{sec:sn}

Corollary \ref{bornes} shows the compactness of the set $\{\Omega\in{\cal C}_{V_0}/\,\lambda_1(\Omega)\leq M\}$. This implies that $\lambda_1$ is coercive. Actually, ${\cal C}_{V_0}$ itself is compact (see section \ref{sec:Bengurian}).

\subsubsection{Case $\delta{=}0$.}
\label{sec:rn}

In this case, $\{\Omega\in{\cal C}_{V_0}/\,\lambda_1(\Omega)\leq M\}$ is also  a compact set and so $\lambda_1$ is coercive. Indeed by Proposition \ref{inrad}, a convex domain $\Omega$ in this set contains the ball $B(x_0,\frac{\pi}{2\sqrt{M+n-1}})$. Set $y\in\Omega$ such that $d(x_0,y)={\rm diam}\,\Omega/4$. Since $\Omega$ is convex, it contains the convex hull of $B(x_0,\frac{\pi}{2\sqrt{M}})\cup\{y\}$ whose volume must be bounded from above by $V_0$. We deduce that ${\rm diam}\, \Omega$ is bounded from above by a function of $M$ and $V_0$. We conclude by Corollary \ref{bornes}.

\subsubsection{Case $\delta{=}-1$.}
\label{sec:hn}

We cannot argue as easily as in the previous case because in $\mathbb{H}^n$, the volume of the convex hull of $B(x_0,\frac{\pi}{2\sqrt{M+n-1}})\cup\{y\}$ does not tend to $\infty$ with $d(x_0,y)$. 
We will prove simultaneously the coercivity of $\lambda_1$ and the property 
\begin{equation}\label{truc}
 \inf_{\mathcal{C'}}\lambda_1(\om) > \lambda_1^*(V_0),
 \end{equation} 
where $\mathcal{C'}=\{\om$ unbounded convex sets; $vol(\om)=V_0\}$. 

These two facts prove Theorem \ref{stabFK}. First, we need to establish some lemmata.

\begin{lemma}Let $\om$ be a domain of a complete Riemannian manifold $(M^n,g)$. Then for any $R\geq1$, $\alpha,\,\gamma\in]0,1[$ and $y_0\in M$, we have
$$\displaylines{\min\Bigl[\lambda_1\bigl(\Omega\cap B(y_0,R)\bigr),\lambda_1\bigl(\Omega\setminus B(y_0,\gamma R)\bigr)\Bigr]\hfill\cr\hfill\leq\frac{1}{(1-R^{-\alpha})^2}\Bigl[\lambda_1(\Omega)+\frac{8}{(1-\gamma)^2R^{2(1-\alpha)}}\Bigr]}$$
where $\lambda_1$ stands for the bottom of the spectrum, $\Omega$ can be of infinite volume and we have set $\lambda_1(\emptyset)=\infty$.
\end{lemma}

\begin{proof}
The proof relies on the variational characterization of the first eigenvalue.

We set $N=E(R^\alpha)+1$, $B_r=B(y_0,r)$, $A_{r,r'}=\Omega\cap(B_r\setminus B_{r'})$ and $r_k=\gamma R+(1-\gamma)R\frac{k}{N}$ for any integer $0\leq k\leq N$. Then, for any $u\in H^1_0(\Omega)$, we have
$$
  \int_{\Omega} u^2\geq\sum_{k=0}^{N-1}\int_{A_{r_{k+1},\,r_k}}^{~}u^2\geq N\int_{A_{r_{k_0+1},\,r_{k_0}}}^{~}u^2
$$
for at least one integer $0\leq k_0\leq N-1$. Let $\phi$ and $\psi$ be the two functions defined on $\R^+$ by:

--$\phi$ is non-decreasing, $\phi=0$ on $[0,\frac{r_{k_0}+r_{k_0+1}}{2}]$, $\phi=1$ on $[r_{k_0+1},\infty[$ and $\|\nabla \phi\|_\infty\leq\frac{2N}{(1-\gamma)R}$,

--$\psi$ is non-increasing, $\psi=1$ on $[0,r_{k_0}]$, $\psi=0$ on $[\frac{r_{k_0}+r_{k_0+1}}{2},\infty[$ and $\|\nabla \psi\|_\infty\leq\frac{2N}{(1-\gamma)R}$,

For $g(x)=\psi\bigl(d(y_0,x)\bigr) u(x)$ and $h(x)=\phi\bigl(d(y_0,x)\bigr) u(x)$, we have
$$\displaylines{
  \int_{\Omega\cap B_{\frac{r_{k_0}+r_{k_0+1}}{2}}} g^2+\int_{\Omega\setminus B_{\frac{r_{k_0}+r_{k_0+1}}{2}}} h^2=\int_{\Omega}(g+h)^2\cr
\geq\int_{\Omega}u^2-\int_{A_{r_{k_0+1},\,r_{k_0}}}^{~}u^2\geq\frac{N-1}{N}\int_{\Omega}u^2}
$$
Since $|dg+dh|^2=|(\psi+\phi)du+ud(\psi+\phi)|^2$, we obtain
$$\displaylines{
  \int_{\Omega\cap B_{\frac{r_{k_0}+r_{k_0+1}}{2}}} |dg|^2+\int_{\Omega\setminus B_{\frac{r_{k_0}+r_{k_0+1}}{2}}} |dh|^2=\int_{\Omega}|dg+dh|^2\hfill\cr
\leq(1+R^{-\alpha})\int_{\Omega}(\phi+\psi)^2|du|^2+(1+R^\alpha)\int_{\Omega}^{~}u^2|d\phi+d\psi|^2\cr
\leq(1+R^{-\alpha})\int_{\Omega}|du|^2+(1+R^\alpha)\int_{A_{r_{k_0+1},\,r_{k_0}}}^{~}u^2\frac{4N^2}{(1-\gamma)^2R^2}\cr
\leq(1+R^{-\alpha})\int_{\Omega}|du|^2+(1+R^\alpha)\int_{\Omega}^{~}u^2\frac{4N}{(1-\gamma)^2R^2}
.}$$
We infer
$$\displaylines{\min\Bigl[\lambda_1\bigl(\Omega\cap B(y_0,R)\bigr),\lambda_1(\Omega\setminus B(y_0,\gamma R)\bigr)\Bigr]\hfill\cr
\hfill\leq\min\Bigl(\frac{\int_{\Omega\cap B_{\frac{r_{k_0}+r_{k_0+1}}{2}}} |dg|^2}{\int_{\Omega\cap B_{\frac{r_{k_0}+r_{k_0+1}}{2}}}g^2},\frac{\int_{\Omega\setminus B_{\frac{r_{k_0}+r_{k_0+1}}{2}}} |dh|^2}{\int_{\Omega\setminus B_{\frac{r_{k_0}+r_{k_0+1}}{2}}}h^2}\Bigr)\hfill\cr
\leq\frac{\int_{\Omega\cap B_{\frac{r_{k_0}+r_{k_0+1}}{2}}} |dg|^2+\int_{\Omega\setminus B_{\frac{r_{k_0}+r_{k_0+1}}{2}}} |dh|^2}{\int_{\Omega\cap B_{\frac{r_{k_0}+r_{k_0+1}}{2}}} g^2+\int_{\Omega\setminus B_{\frac{r_{k_0}+r_{k_0+1}}{2}}}h^2}\cr
\leq\frac{1}{(1-R^{-\alpha})^2}\Bigl[\frac{\int_{\Omega}|du|^2}{\int_{\Omega}u^2}+\frac{8}{(1-\gamma)^2R^{2(1-\alpha)}}\Bigr]}$$
\end{proof}

This lemma implies the following result.
\begin{lemma}
  For any $V_0>0$ there exist $C(V_0,n)>\lambda_1^*(V_0)$ and $R(V_0,n)>0$ such that, for any bounded convex set $\Omega$ which satisfies $vol (\om) \leq V_0$ and $\lambda_1(\Omega)\in[\lambda_1^*(V_0),C(V_0,n)]$, we have
$$\lambda_1(\Omega)\leq\lambda_1\bigl(\Omega\cap B(x_0,R)\bigr)\leq\frac{1}{(1-R^{-\alpha})^2}\Bigl[\lambda_1(\Omega)+\frac{32}{R^{2(1-\alpha)}}\Bigr]$$
for any $\alpha\in]0,1[$, $R\geq R(V_0,n,\alpha)$ and $x_0$ such that $B\bigr(x_0,{\rm Inrad(\Omega)}\bigl)\subset \overline{\Omega}$.
\end{lemma}

\begin{proof}
  We set $r(V_0,n)=\frac{\pi}{\sqrt{2\lambda_1^*(V_0)+n-1}}$ and 
$$C(V_0,n)=\min\bigl[2\lambda_1^*(V_0),\frac{\lambda_1^*(V_0)+\lambda_1^*\bigl(V_0-\Vol B(x_0,r(V_0,n)/2)\bigr)}{2}\bigr].$$
Then, by Proposition \ref{inrad}, we have $B(x_0,r(V_0,n)/2)\subset\Omega$ and so $\Vol\bigl(\Omega\setminus B(x_0,R/2)\bigr)\leq V_0- \Vol B(x_0,r(V_0,n)/2)$ for any $R\geq r(V_0,n)$. By the Faber-Krahn inequality, this implies that $\lambda_1\bigl(\Omega\setminus B(x_0,R/2)\bigr)$ is larger than $\lambda_1^*\bigl(V_0-\Vol B(x_0,r(V_0,n)/2)\bigr)$. Now, we can choose $R(V_0,n)$ large enough in order to have 
$\frac{1}{(1-R^{-\alpha})^2}\Bigl[C(V_0,n)+\frac{32}{R^{2(1-\alpha)}}\Bigr]\leq\lambda_1\bigl(\Omega\setminus B(x_0,R/2)\bigr)$ for any $R\geq R(V_0,n)$. Lemma 5.1 then applies.
\end{proof}
\medskip

Now, we prove (simultaneously) the coercivity of $\lambda_1$ and (\ref{truc}). By definition of the bottom of the spectrum, it is sufficient to prove that every sequence of bounded convex domains $(\om_i)_{i \in \N}$ such that $vol (\om_i) \leq V_0$ and $\lim_i \lambda_1(\om_i)=\lambda_1^*(V_0)$, converges, up to isometries and extraction, to $B(x_0,r_0)$ ($vol (B(x_0,r_0))=V_0$).

Let $(\om_i)_{i \in \N}$ be such a sequence. Up to isometries, we can suppose that the fixed point $x_0\in \hn$ satisfies the condition $B\bigl(x_0,{\rm Inrad}(\Omega_i)\bigr)\subset\Omega_i$ for every $i$. By the lemma above and Corollary \ref{bornes}, the sequence $\bigl(\Omega_i\cap B_R\bigr)_{i \in \mathbb{N}}$ is precompact in ${\cal C}$ for all $R\geq R(V_0,n)$. Up to a diagonal extraction, we can now suppose that for any $n\in\N$ the sequence $\bigl(\Omega_i\cap B_n\bigr)_{i \in \mathbb{N}}$ converges to an element $U_n$ of ${\cal C}$. Using the continuity of $\lambda_1$ and the volume on ${\cal C}$, we get

$$\displaylines{\lambda_1^*(V_0)\leq\lambda_1(U_n)\leq\frac{1}{(1-n^{-1/2})^2}\Bigl[ \lambda_1^*(V_0)+\frac{32}{n}\Bigr]\cr
\mbox{and}\quad\Vol(U_n)\leq V_0}$$
So $\lambda_1(U_n)$ tends to $\lambda_1^*(V_0)$ and by the Faber-Krahn inequality, we must have $\Vol U_n\to V_0$.
Moreover, $(U_n)_{n \in \mathbb{N}}$ is a nondecreasing sequence of convex sets for the inclusion. As a consequence,  $\Omega=\bigcup_n U_n$ is a convex domain of volume $V_0$ and first eigenvalue $\lambda_1(\om)=\lambda_1^*(V_0)$. By Proposition \ref{FK}, $\Omega=B(x_0,r_0)$ and we infer that the sequence $(\Omega_i)_{i \in \N}$ converges to $B(x_0,r_0)$ in ${\cal C}$.

\section{Coercivity of the $\lambda_1/\lambda_2$ functional}
\label{sec:Bengurian}
\subsection{case $\delta=0$}
We show that, on ${\cal C}_{V_0}$, $\lambda_1/\lambda_2$ tends to $1$ when $\lambda_1$ tends to $\infty$. By section \ref{sec:rn},  $\inf_{{\cal C}_{V_0}}\lambda_1/\lambda_2<1$ and the fact that $\lambda_1/\lambda_2$ is invariant under homothetie on the domains, this implies Theorem \ref{stabB} in $\rn$.

By a classical result due to Jones, for any $\Omega\in{\cal C}_{V_0}$ there exists an ellipsoid ${\cal E}$ such that ${\cal E}\subset\Omega\subset\sqrt{n}{\cal E}$. We easily infer that there is a n-rectangle $R$ with edges of lengths $L_1\leq\cdots\leq L_n$,  such that $R\subset\Omega\subset nR$. This gives

\begin{equation}\label{vol}
\lambda_1(\Omega)\leq\lambda_1(R)\leq\frac{n\pi^2}{L_1^2}\quad\mbox{ and }\quad V_0\leq L_n^{n-1}n^nL_1\text,
\end{equation}
and so $L_n\geq\Bigl(\frac{V_0\sqrt{\lambda_1}}{\pi n^{n+\frac{1}{2}}}\Bigr)^\frac{1}{n-1}$.
Following \cite{Jer}, we can translate and rotate $\Omega$ so that $R$ be centred in $(0,\ldots,0)$ and the edge of $R$ of length $L_n$ be parallel to the last coordinate axis. We denote $\Omega(y)=\{x\in\R^{n-1}/\,(x,y)\in\Omega\}$ and $\lambda(y)=\lambda_1\bigl(\Omega(y)\bigr)$. Then, if $f$ is an eigenfunction of $\Omega$ associated to the first eigenvalue, we have
$$\displaylines{\hfill\int_\Omega f^2=\int_\Omega\frac{|\nabla f|^2}{\lambda_1(\Omega)}\geq\int_\R\int_{\Omega(y)}\frac{|\nabla_x f(x,y)|^2}{\lambda_1(\Omega)}dxdy\hfill\cr
\hfill\geq\int_\R\frac{\lambda_1\bigl(\Omega(y)\bigr)}{\lambda_1(\Omega)}\int_{\Omega(y)}|f(x,y)|^2dxdy\hfill}$$
Thus, there is $y$ such that $\lambda_1(\Omega)\geq\lambda_1\bigl(\Omega(y)\bigr)$. By convexity of $\Omega$ we deduce that $\bigl(1-(\frac{L_n}{2})^{-\frac{2}{3}}\bigr)\Omega(y)\times[y-(\frac{L_n}{2})^\frac{1}{3},y+(\frac{L_n}{2})^\frac{1}{3}]$ is contained in $\Omega$ and consequently,
$$\lambda_1(\Omega)\leq\lambda_2(\Omega)\leq\lambda_2\Bigl(\bigl(1-(\frac{L_n}{2})^{-\frac{2}{3}}\bigr)\Omega(y)\times[y-(\frac{L_n}{2})^\frac{1}{3},y+(\frac{L_n}{2})^\frac{1}{3}]\Bigr)$$
\begin{equation}\label{norml}
 \leq\frac{\lambda_1\bigl(\Omega(y)\bigr)}{\bigl(1-(\frac{L_n}{2})^{-\frac{2}{3}}\bigr)^2}+\frac{2\pi^2}{(\frac{L_n}{2})^\frac{2}{3}}\leq\frac{\lambda_1(\Omega)}{\bigl(1-(\frac{L_n}{2})^{-\frac{2}{3}}\bigr)^2}+\frac{2\pi^2}{(\frac{L_n}{2})^\frac{2}{3}}.
\end{equation}

Since we have shown above that $L_n\to\infty$ when $\lambda_1\to\infty$, we obtain that $\lambda_1/\lambda_2$ tend to $1$ when $\lambda_1$ tends to $\infty$.

\begin{remark}
  The same method could be used to show that for any integers $p\leq q$, $\lambda_p/\lambda_q$ tends to $1$ when $\lambda_1$ tends to $\infty$ on ${\cal C}_{V_0}$. We conlude that for any $p\leq q$ there exists a convex domain (to be determined) which minimizes the quotient $\lambda_p/\lambda_q$.
\end{remark}

\begin{remark}
The inequations \eqref{norml} imply that for any convex domain $\Omega$ of $\rn$, $x_0\in\Omega$ and $R>0$ such that $B(x_0,R)$ does not contain $\Omega$, we have 
$$\lambda_1(\Omega\cap B(x_0,R))\leq\lambda_2(\Omega\cap B(x_0,R))\leq\frac{\lambda_1(\Omega\cap B(x_0,R))}{\bigl(1-C(n)R^{-\frac{2}{3}}\bigr)^2}+\frac{C(n)}{R^\frac{2}{3}}$$
and so $\lambda_2(\Omega\cap B(x_0,R))$ tends to $\lambda_1(\Omega)$ when $R$ tends to $\infty$. 
We conclude that a convex Euclidean domain with a spectral gap is bounded (hence has a discrete spectrum) and that its diameter is bounded from above by
$C(n)\bigl(\frac{1+\lambda_1}{\lambda_2-\lambda_1}\bigr)^{3/2}$. This implies readily the coercivity of $\lambda_1/\lambda_2$ on the set of convex Euclidean domains of fixed $\lambda_1$, from which we infer Theorem \ref{stabBL} in $\rn$.
\end{remark}

\subsection{Case $\delta=1$}
The coercivity $\lambda_1/\lambda_2$ on the set of convex domains with $\lambda_1=\lambda$ follows from Lemma \ref{bornes}. On ${\cal C}_{V_0}$, this follows from the compactness of ${\cal C}_{V_0}$ which, by Corollary \ref{bornes}, is a consequence of the inequality ${\rm Inrad}\,(\Omega)\geq C(n)\Vol(\Omega)$. This inequality holds true for any convex domain of $\sn$ as explained below.

First, using the inequality (\ref{vol}), based on the Jones ellipsoid, we get easily that for any convex domain contained in a geodesic ball of radius $R$ in $\R^n$, we have $\Vol(\Omega)\leq n^n R^{n-1}{\rm Inrad}\, (\Omega)$. Now, since $\S^n$ can be covered by $2(n+1)$ balls of radius $R_n={\rm arccos}(\frac{1}{\sqrt{n+1}})$ we infer that there is a point $x_0$ in $\S^n$ such that $\Vol \bigl(\Omega\cap B(x_0,R_n)\bigr)\geq\frac{1}{2(n+1)}\Vol(\Omega)$. Using the canonical embedding of $\S^n$ in $\mathbb{R}^{n+1}$, we can project $B(x_0,R_n)$ onto the tangent space $T_{x_0}\S^n$  (using the origin of Euclidean space). This map $P_0$ is a quasi-isometry from the ball $B(x_0,R_n)$ in $\S^n$ to the geodesic ball $B(x_0,\sqrt{n})$ in Euclidean space, which preserves the convexity. Then, we have
$$\displaylines{{\rm Inrad}_{\S^n}\, (\Omega)\geq{\rm Inrad}_{\S^n}\, \bigl(\Omega\cap B(x_0,R_n)\bigr)\hfill\cr\hfill\geq C_1(n){\rm Inrad}_{T_{x_0}\S^n}\, \bigl(P_0\bigl(\Omega\cap B(x_0,R_n)\bigr)\bigr)\hfill\cr\hfill\geq C_2(n)\Vol_{T_{x_0}\S^n} \bigl(P_0\bigl(\Omega\cap B(x_0,R_n)\bigr)\bigr)\geq C(n)\Vol_{\S^n}(\Omega)}$$

\subsection{Case $\delta=-1$.}

In this section, we prove simultaneously the coercivity of the functional $\lambda_1/\lambda_2$ on bounded convex domains whose first eigenvalue is fixed and the property
\begin{equation}\label{truc2}
 \sup_{\mathcal{C'}}\lambda_2(\om) < \lambda_2^*(\lambda),
 \end{equation}
where $\mathcal{C'}=\{\om$ unbounded convex sets; $\lambda_1(\om)=\lambda\}$.

\medskip 

\noindent These two properties imply Theorem \ref{stabBL}.

\medskip

We need the following result whose proof follows easily from the min-max principle (see \cite[theorem XIII.1]{RS}).
\begin{lemma}\label{le}
 Let $\Omega$ be a convex domain in $\mathbb{H}^n$ such that the bottom of the spectrum is an eigenvalue. 
 Then for any fixed point $x_0\in\hn$, we have
$$\lim_{R\to\infty}\lambda_i(\Omega\cap B(x_0,R))=\lambda_i(\Omega),\quad\mbox{for }i=1,2\text.$$
\end{lemma}
Thanks to this lemma, the coercivity property  and  the inequality (\ref{truc2}) reduce to the fact below.

\smallskip

\noindent Every sequence $(\om_i)_{i\in \N} \in \mathcal{C}$ such that $\lim_i \lambda_1(\om_i)=\lambda$ and $\lim_i \lambda_2(\om_i) \\ =\lambda_2^*(\lambda)$, converges (up to extraction) to a ball such that $\lambda_1(B)=\lambda$.




\bigskip

First, we show that a lower bound on the spectral gap implies some estimates on the first eigenfunction.

\begin{lemma}
Let $\Omega$ be a bounded domain of $\hn$. If $u\in H^1_0(\Omega)$ satisfies $\lap u=\lambda_1(\Omega)u $ then there is a point $x_m\in\hn$ such that
$$\Bigl(\lambda_2(\Omega)-\lambda_1(\Omega)-\frac{n{-}1}{\sinh^2(R)}\Bigr)\int_{\Omega\setminus B(x_m,R)}u^2\leq \frac{n}{R^2}\int_{\Omega\cap B(x_m,R)}u^2$$
for any $R>0$. This implies, for any $R\geq2\sqrt{\frac{n-1}{\lambda_2(\Omega)-\lambda_1(\Omega)}}$,
$$\displaylines{\lambda_1(\Omega)\leq\lambda_1\bigl(\Omega\cap B(x_m,R)\bigr)\hfill\cr
\hfill\leq\frac{(1+\frac{1}{R^2})}{1-\frac{4n}{(\lambda_2(\Omega)-\lambda_1(\Omega))R^2+4})}\Bigl(\lambda_1(\Omega)+\frac{4n}{(\lambda_2(\Omega)-\lambda_1(\Omega))R^2+4}\Bigr).}$$
\end{lemma}

\begin{proof}
  Proposition \ref{specgap} applied to $g(s)=s/R$ on $[0,R]$ and $g(s)=1$ on $[R,\infty[$ gives a point $x_m\in\hn$ such that
$$\displaylines{\lambda_2(\Omega)-\lambda_1(\Omega)
\leq\frac{\frac{n}{R^2}\int_{\Omega\cap B(x_m,R)}u^2(x)dx+\frac{n-1}{\sinh^2R}\int_{\Omega\setminus B(x_m,R)}u^2(x)dx}{\int_{\Omega\setminus B(x_m,R)} u^2(x)dx}}$$
which gives the first estimate. 

For the second estimate, we set
$\psi$ the non-increasing Lipschitzian function defined on $\R^+$ by $\psi=1$ on 
$[0,R/2]$, $\psi=0$ on $[R,\infty[$ and $\|\nabla\psi\|\infty=\frac{2}{R}$.
Then, we have 
\begin{eqnarray*}
 |d(\psi u)|^2&=&\psi^2|du|^2+2u\psi(d\psi,du)+u^2|d\psi|^2\cr
&\leq&(1+\frac{1}{R^2})\psi^2|du|^2+(1+R^2)|d\psi|^2u^2
\end{eqnarray*}
So, we infer

$$\displaylines{\lambda_1\bigl(\Omega\cap B(x_m,R)\bigr)\leq\frac{\int_{\Omega\cap B(x_m,R)}|d(\psi u)|^2}{\int_{\Omega\cap B(x_m,R)}(\psi u)^2}\hfill\cr
\hfill\leq(1+\frac{1}{R^2})\frac{\int_{\Omega}\psi^2|du|^2}{\int_{\Omega\cap B(x_m,R/2)}u^2}+(1+R^2)\frac{\int_{\Omega}|d\psi|^2u^2}{\int_{\Omega\cap B(x_m,R/2)}u^2}\cr
\leq(1+\frac{1}{R^2})\frac{\int_{\Omega}|du|^2}{\int_{\Omega}u^2}\Bigl(1+\frac{\int_{\Omega\setminus B(x_m,R/2)}u^2}{\int_{\Omega\cap B(x_m,R/2)}u^2}\Bigr)+4(1+\frac{1}{R^2})\frac{\int_{\Omega\setminus B(x_m,R/2)}u^2}{\int_{\Omega\cap B(x_m,R/2)}u^2}.}$$
By the first estimate, we have
$$\frac{\int_{\Omega\setminus B(x_m,R/2)}u^2}{\int_{\Omega\cap B(x_m,R/2)}u^2}\leq\frac{4n}{(\lambda_2(\Omega)-\lambda_1(\Omega))R^2-4(n-1)}$$
from which we infer
$$\displaylines{\lambda_1\bigl(\Omega\cap B(x_m,R)\bigr)\leq\hfill\cr
\hfill\frac{(1+\frac{1}{R^2})\bigl((\lambda_2(\Omega)-\lambda_1(\Omega))R^2+4\bigr)}{(\lambda_2(\Omega)-\lambda_1(\Omega))R^2-4(n-1)}\Bigl(\lambda_1(\Omega)+\frac{4n}{(\lambda_2(\Omega)-\lambda_1(\Omega))R^2+4}\Bigr)}$$
\end{proof}

\medskip

Let $(\Omega_i)_{i \in \N}\in{\cal C}$ such that $\lim_i \lambda_1(\Omega_i)=\lambda$ and $\lim_i\lambda_2(\Omega_i)=\lambda_2^*(\lambda)$. We can assume that $\lambda_2(\Omega_i)-\lambda_1(\Omega_i)>\frac{\lambda_2^*(\lambda)-\lambda}{2}>0$. Note that by the preceeding lemma, we infer that for any $R\geq4\sqrt{\frac{n-1}{\lambda_2^*(\lambda)-\lambda}}$ we have $\lambda_1\bigl(\Omega_i\cap B(x_m^i,R)\bigr)\leq C(\lambda,n)$ (where $C(\lambda,n)$ is a universal funtion and $x_m^i$ is the center of mass of $\Omega_i$). This implies, by Proposition \ref{inrad}, that we can suppose (up to isometry) $x_m^i\in B\bigl(x_0,4\sqrt{\frac{n-1}{\lambda_2^*(\lambda)-\lambda}}\bigr)$ and $B\bigl(x_0,r(\lambda,n)\bigr)\subset\Omega_i$ 
for all $i$. Then, the sequence $\bigl(\Omega_i\cap B(x_0,R)\bigr)$ is included in a compact set of ${\cal C}$ (see Corollary \ref{bornes}). By diagonal extraction, we can suppose that for any $k\in\N$ the sequence $\bigl(\Omega_i\cap B(x_0,k)\bigr)_{i \in \N}$ converges to an element $U_k$ of ${\cal C}$. By continuity of $\lambda_1$ on ${\cal C}$, we have
$$\displaylines{\lambda\leq\lambda_1(U_k)=\lim_i\lambda_1\bigl(\Omega_i\cap B(x_0,k)\bigr)\hfill\cr
\hfill\leq\lim_i\lambda_1\bigl(\Omega_i\cap
 B(x_m^i,k-4\sqrt{\frac{n-1}{\lambda_2^*(\lambda)-\lambda}})\bigr)\leq f(k,\lambda,n)}$$
where $f(k,\lambda,n)$ is a universal function given by the preceding lemma and that converge to $\lambda$ when $k$ tends to $\infty$.
So $\lambda_1(U_k)$ tends to $\lambda$. As in the subsection \ref{sec:hn}, we set $\Omega=\cup_kU_k$. Then $\Omega$ is a convex domain with $\lambda_1(\Omega)=\lim_k\lambda_1(U_k)=\lambda \,$
(since $U_k=\om \cap B(x_0,k)$) and 
$$\lambda_2(\Omega)=\lim_k\lambda_2(U_k)=\lim_k\lim_i\lambda_2(\Omega_i\cap B(x_0,k))\geq\lim_i\lambda_2(\Omega_i)=\lambda_2^*(\lambda).$$
Then, we conclude by Theorem \ref{conv2}.
\section*{Appendices}
\begin{appendix}
  
\section{Trigonometric computations}
\label{supp}

In this appendix, we perform the computations of $l'(0)$ and $L'(0)$ used in the proof of Lemma \ref{support}. We denote by $\beta$ the half angle at $y_0$ of the geodesic double cone tangent to the ball $B(x_0,r)$. By the law of sines, we have $\sin\beta=\frac{s_\delta(r)}{s_\delta\bigl(d(x_0,y_0)\bigr)}$ and $\frac{s_\delta\bigl(l(t)\bigr)}{\sin\beta}=\frac{s_\delta\bigl(l_1(t)\bigr)}{\sin t}$, where we have set $l_1\bigl(d(u_0,v)\bigr)=d(x_0,Z(v))$. By letting $t$ tend to $0$, we get $l_1'(0)=\frac{s_\delta^2\bigl(d(x_0,y_0)\bigr)}{s_\delta(r)}$. On the other hand, the cosine law gives us the equation $c_\delta(l)=c_\delta(l_1)c_\delta\bigl(d(x_0,y_0)\bigr)+\delta s_\delta(l_1)s_\delta\bigl(d(x_0,y_0)\bigr)\cos\beta$ (resp. $l^2=l_1^2+\bigl(d(x_0,y_0)\bigr)^2-2l_1d(x_0,y_0)\cos\beta$ if $\delta=0$), whose derivative at $t=0$ gives the relation $l'(0)=-l_1'(0)\cos\beta$. We easily deduce the relation $l'(0)=-s_\delta\bigl(d(x_0,y_0)\bigr)\sqrt{\Bigl(\frac{s_\delta\bigl(d(x_0,y_0)\bigr)}{s_\delta(r)}\Bigr)^2-1}$. Note that for $L'(0)$ we just have to replace $\beta$ by $\pi-\beta$ in what preceed.

\section{Moser's iteration}
\label{iter}

In this section, we prove the inequality \eqref{itera} used in the proof of lemma \ref{KOMPAQ}.

 Set  $0\leq\eta\leq1$ a ${\cal C}^1$ function such that $\eta\equiv 1$ on $B(x_0,\alpha r)$ (for $\alpha\in]0,1[$ and $1\geq r>0$), $\eta\equiv 0$ on $X^\delta\setminus B(x_0,r)$ and $|d\eta|\leq2/(1-\alpha)r$. 

We fix $m>0$ and $\beta\geq0$ and set $h=\inf(m,f_1)$, $u=f_1h^\frac{\beta}{2}$ and $\phi=\eta^2h^\beta f_1\in H_0^1(\Omega)$. Then we have
$$\displaylines{\lambda_1\int_\Omega\eta^2u^2=\lambda_1\int_\Omega\phi f_1\geq\int_\Omega(df_1,d\phi)\hfill\cr
\hfill\geq\beta\int_\Omega\eta^2h^\beta |dh|^2+\frac{1}{2}\int_\Omega\eta^2h^\beta|df_1|^2-2\int_\Omega|d\eta|^2h^\beta f_1^2,}$$
where we used $2\eta f_1(df_1,d\eta)\geq-\frac{1}{2}\eta^2|df_1|^2-2f_1^2|dh|^2$.
This inequality, combined with the inequalities $|d(\eta u)|^2\leq2u^2|d\eta|^2+2\eta^2|du|^2$ and $|du|^2\leq(1+\beta)h^\beta(2\beta|dh|^2+|df_1|^2)$, gives
$$\int_\Omega|d(u\eta)|^2\leq(10+4\lambda_1)(1+\beta)\int_\Omega u^2(\eta^2+|d\eta|^2).$$
Hence the Sobolev inequality $(*)$ applied to $u\eta$ implies
$$\Bigl(\int_{B(x_0,\alpha r)}h^\frac{(2+\beta)n}{n-2}\Bigr)^\frac{n-2}{n}\leq \frac{5C(n)(10+4\lambda_1)(1+\beta)}{(1-\alpha)^2r^2}\int_{B(x_0,r)}f_1^{(2+\beta)}$$
Then, we let $m$ tend to $\infty$ and set $r_k=\frac{1}{2^{\sqrt{k}}}$, $\alpha_k=2^{\sqrt{k}-\sqrt{k+1}}$ and $\beta_k=2\bigl(\frac{n}{n-2}\bigr)^k-2$. By multiplying the $(2+\beta_k)$-th square root of the inequalities obtained for $1\leq k\leq K-1$ we infer
$$\Bigl(\int_{B(x_0,r_K)}h^{2(\frac{n}{n-2})^K}\Bigr)^\frac{1}{2(n/n-2)^K}\leq A(n,K)(1+\lambda_1)^{\gamma(K)}\int_{B(x_0,r)}f_1.$$
By definition of $r_K$, we have $\Bigl[\int_{B(x_0,r_K)}f^{2(\frac{n}{n-2})^K}\Bigr]^{\frac{1}{2(n/n-2)^K}}$ tends to $f(x_0)$ when $K$ tends to $+\infty$, meanwhile $A(n,K)$ and $\gamma(K)$ converge, which gives \eqref{itera}.

\section{Proof of lemma \ref{masscenter}}
\label{proofcenter}

  This lemma is essentially proven in \cite{ash1,ash2,BL} for $u$ with compact support (which includes the case $\delta=1$) but we need to apply it to an eigenfunction $u$ of a convex (unbounded) domain $\Omega$ ,and so we have to extend it in the case $\delta=0,-1$.

In the sequel of the proof, $X$ denotes $\rn$ or $\hn$. We fix $x_0\in X$ and define
\begin{eqnarray*}
  F~:~T_{x_0}X&\to&T_{x_0}X\\
v&\mapsto&d({\rm exp}_{x_0}^{-1})\Bigl(\int_Xg\bigl(d(\bar{v},y)\bigr)\frac{{\rm exp}_{\bar{v}}^{-1}(y)}{d(\bar{v},y)}u^2(y)\,dy\Bigr)
\end{eqnarray*}
where we have set $\bar{v}={\rm exp}_{x_0}(v)$. We set $m=\liminf_{+\infty} g$. Let $R_1>0$ such that $\int_{X\setminus B(x_0,R_1)}u^2\leq\min\bigl(\frac{m}{32\|g\|_\infty},\frac{1}{2}\bigr)$. Then for any $v\in T_{x_0}X$ with $|v|\geq R_1$ we easily have
$$\Bigl|F(v)-d({\rm exp}_{x_0}^{-1})\Bigl(\int_{B(x_0,R_1)}g\bigl(d(\bar{v},y)\bigr)\frac{{\rm exp}_{\bar{v}}^{-1}(y)}{d(\bar{v},y)}u^2(y)\,dy\Bigr)\Bigr|\leq\frac{m}{32}.$$
Note that $d({\rm exp}_{x_0}^{-1})\circ{\rm exp}_{\bar{v}}^{-1}(x_0)=-v$ and so we infer that for any $v\in T_{x_0}X$ with $|v|\geq R_1$ we have
$$\Bigl|F(v)+\lambda(v)\frac{v}{|v|}\Bigr|\leq\|g\|_\infty\int_{B(x_0,R_1)}\bigl|\frac{{\rm exp}_{\bar{v}}^{-1}(x_0)}{d(\bar{v},x_0)}-\frac{{\rm exp}_{\bar{v}}^{-1}(y)}{d(\bar{v},y)}\bigr|\,dy+\frac{m}{32}$$
where we have set $\lambda(v)=\int_{B(x_0,R_1)}g\bigl(d(\bar{v},y)\bigr)u^2(y)\,dy$ and used the fact that $d({\rm exp}_{x_0}^{-1})$ is a contraction. Then we have $\lambda(v)\geq\frac{m}{4}>0$ for any $v$ with $|v|\geq R_2\geq R_1$. Note also that the integrand above measures the difference between the unit tangent vectors at $\bar{v}$ to the minimizing geodesic from $\bar{v}$ to $x_0$ and $y\in B(x_0,R_1)$. By the law of cosines, we can easily show that this quantity uniformly tends to zero on $B(x_0,R_1)$ when $|v|$ tends to $+\infty$. Hence, there exists $R_3>0$ such that  for any $v\in T_{x_0}X$ which satisfies $|v|\geq R_3$, we have
$$\Bigl|F(v)+\lambda(v)\frac{v}{|v|}\Bigr|\leq\frac{m}{16}\quad\mbox{and}\quad\lambda(v)\geq\frac{m}{4}$$
We have to show that $F$ is zero somewhere. If not, the following map (with $R>R_3$) is well-defined.
\begin{eqnarray*}
  G~:~B(0,1)\subset T_{x_0}X&\to&S(0,1)\subset T_{x_0}X\\
v&\mapsto&\frac{F(-Rv)}{|F(-Rv)|}.
\end{eqnarray*}
Moreover, the map $G$ is continuous and satisfies \\
$|G(v)-v|\leq\frac{2\bigl|F(-Rv)+\lambda(-Rv)\bigr|}{|F(-Rv)|}\leq 4/3$ for any $v\in S(0,1)$. So, we could then easily construct a retraction from $B(0,2)$ to $S(0,2)$.

\section{A result of Li and Yau}
\label{LY}
\begin{lemma}[Li-Yau]
  Let $\Omega$ be a bounded and smooth domain with positive mean curvature (for the exterior normal). If $f$ an eigenfunction associated to the first eigenvalue of the Dirichlet problem on $\Omega$, then we have,
$$|\nabla f|^2\leq\lambda_1(\|f\|_\infty^2-f^2)$$
(resp.
$$|\nabla f|^2\leq(\lambda_1+n-1)(\|f\|_\infty^2-f^2)$$ if $\delta=-1$).
\end{lemma}
\begin{proof}
  Let $F=\frac{|\nabla f|^2}{\beta-f^2}$ where $\beta=(1+\epsilon)\|f\|_\infty^2$. Then we have
$$dF(v)=\frac{2|\nabla f|^2}{\beta-f^2}\Bigl(\frac{\Hess f(\frac{\nabla f}{|\nabla f|},v)}{|\nabla f|}+\frac{fdf(v)}{\beta-f^2}\Bigr).$$
If $x_0$ is a point of $\partial \Omega$ then $\nu=\nabla f(x_0)/|\nabla f(x_0)|$ is well-defined (by the strong maximum principle applied to $f$) and is the interior normal to $\Omega$ at $x_0$. We then have
$$dF(\nu)=\frac{2|\nabla f|^2}{\beta-f^2}\Bigl(\frac{\Hess f (\nu,\nu)}{|\nabla f|}+\frac{f|\nabla f|}{\beta-f^2}\Bigr)\geq0,$$
since $\frac{\Hess f (\nu,\nu)}{|\nabla f|}=-\frac{\lap f}{|\nabla f|}+\mu(x_0)$, where $\mu(x_0)$ is the mean curvature of $\partial \Omega$ at $x_0$. We infer by the strong maximum principle that at a point $x_0$ where $F$ reaches its maximum on $\bar{\Omega}$ we must have 
$$\displaylines{\hfill dF(x_0)=0\hfill\mbox{and}\hfill \lap F (x_0)\geq0\hfill}$$
The first equation and our computation of $dF$ imply that $\nabla f/|\nabla f|$ is an eigenvector of $\Hess f(x_0)$ with respect to $g(x_0)$, associated to the eigenvalue $-\frac{f|\nabla f|^2}{\beta-f^2}$ (we have $\nabla f(x_0)\neq 0$ since $F\neq 0$). So we have $|\Hess f(x_0)|^2\geq f^2F^2$.

From the Bochner formula $\frac{1}{2}\lap |\nabla f|^2=\lambda_1|\nabla f|^2-|\Hess f|^2-\Ric(\nabla f,\nabla f)$ (where $\Ric$ denote the Ricci curvature tensor of $X^\delta$) we infer that, at $x_0$, we have
$$\displaylines{|\nabla f|^2F-f\lap f F+\frac{(\beta-f^2)}{2}\lap F=\frac{1}{2}\lap \bigl((\beta-f^2)F\bigr)\hfill\cr
\hfill\leq\lambda_1|\nabla f|^2-f^2F^2-\delta(n-1)|\nabla f|^2.}$$
Since $\lap F(x_0)\geq0$ and $|\nabla f|^2=F(\beta-f^2)$ we readily obtain the estimate $F(x_0)\leq\lambda_1$ (resp. $F(x_0)\leq\lambda_1+n-1$ if $\delta=-1$). Then, we just have to let $\epsilon$ tend to $0$.
\end{proof}

To get Proposition \ref{inrad} in the case of a smooth convex domain, let $f$ denote a positive eigenfunction associated to $\lambda_1$ and $z_0\in\Omega$ a point where $f(z_0)=\|f\|_\infty$. Set $\gamma$ a normal geodesic from $z_0$ to a point $y\in\partial\Omega$. By lemma \ref{LY}, we obtain
$$\bigl(\arcsin( f\circ\gamma/\|f\|_\infty)\bigr)'\geq-\sqrt{\lambda_1+\delta(n-1)}$$
and so that $f\circ\gamma(t)\geq\|f\|_\infty\cos\bigl(\sqrt{\lambda_1+\delta(n-1)}t\bigr)$. Since $f(y)=0$, the geodesic ball $B(z_0,\frac{\pi}{2\sqrt{\lambda_1+\delta(n-1)}})$ is included in $\Omega$.

\end{appendix}

\vfill
\begin{flushleft}
Erwann AUBRY\\
Universit\'{e} de Nice Sophia-Antipolis\\
Laboratoire J.-A. Dieudonn\'{e}\\
UMR6621 (UNSA-CNRS)\\
Parc Valrose\\
F-06108 Nice Cedex (France)\\
eaubry@math.unice.fr
\vspace{0.5cm}

J\'{e}r\^{o}me BERTRAND\\
Scuola Normale Superiore\\
Piazza dei Cavalieri, 7\\
I-56100 Italia \\
j.bertrand@sns.it
\vspace{0.5cm}

Bruno COLBOIS\\
Institut de math\'{e}matiques\\
Universit\'{e} de Neuch\^{a}tel\\
Rue \'Emile Argand, 11\\
Case postale 158\\
CH-2009 Neuch\^{a}tel\\
bruno.colbois@unine.ch

\end{flushleft}

\end{document}